\documentclass[12pt]{article} 
\usepackage{graphicx}
\usepackage{url}

\newcommand{\union}{\cup}
 
\renewcommand{\int}{\cap}

\newcommand{\ci}{\subseteq}

\newcommand{\set}[1]{ {\left\{ {#1} \right\} } } 
 
\newcommand{\seq}[2]{ {{#1}_1 , {#1}_2 , \ldots , {#1}_{#2}} }

\newtheorem{definition}{Definition}
\newtheorem{lemma}[definition]{Lemma}
\newtheorem{theorem}[definition]{Theorem}

\newtheorem{proposition}[definition]{Proposition}
\newtheorem{corollary}[definition]{Corollary}

\newtheorem{observation}[definition]{Observation}

\newcommand{\qed}{\hfill \rule{1ex}{1ex}} 
\newenvironment{proof}{{\bf Proof}: }{\qed}

\newcommand{\C}{\mathbf{C}}

\newcommand{\Z}{\mathbf{Z}}

\newcommand{\ord}{\mathrm{ord}}
\newcommand{\Com}{\mathrm{Com}}
\newcommand{\e}{{\mathrm{e}}}

\title{The enumeration of simple permutations}
\author{
M.H. Albert\thanks{Department of Computer Science,
University of Otago,
Dunedin, New Zealand. 
{\tt malbert@cs.otago.ac.nz}}
\and
M.D. Atkinson\thanks{Department of Computer Science,
University of Otago,
Dunedin, New Zealand. 
{\tt mike@cs.otago.ac.nz}} \and
M. Klazar\thanks{Department of Applied Mathematics (KAM)
and Institute for Theoretical Computer Science (ITI), Charles
University, Malostransk\'e n\'am\v est\'\i\ 25, 118 00 Praha, Czech
Republic. ITI is supported by the project LN00A056 of the Ministry of
Education of the Czech Republic. {\tt
klazar@kam.mff.cuni.cz}}}
\date{}

\begin{document}

\maketitle 

\begin{abstract}
A {\em simple permutation} is one which maps no proper
non-singleton interval onto an interval. We consider
the enumeration of simple permutations from several aspects. Our
results include a straightforward relationship between the ordinary
generating function for simple permutations and that for all
permutations, that the coefficients of this series are not
$P$-recursive, an asymptotic expansion for these coefficients, and a
number of congruence results.

\bigskip

\end{abstract}

\section{Introduction and definitions}

The permutation $2647513$ maps the interval $2..5$ onto the interval
$4..7$.  In other words, it has a \emph{segment} (set of consecutive
positions) whose values form a \emph{range} (set of consecutive
values).  Such a segment is called a \emph{block} of the permutation.
Every permutation has singleton blocks, together with the block
$1..n$.  If these are the only blocks the permutation is called
\emph{simple}.  For example, $58317462$ is simple and the simple
permutations of length up to $5$ are as follows:.

\begin{center}
\begin{tabular}{cl} \hline 
\makebox[2.5cm]{Length} & Simple permutations \\ \hline
\rule{0pt}{16pt}$1$ & 1 \\
$2$ & 12, 21 \\ 
$3$ & None \\ 
$4$ & 2413, 3142 \\ 
$5$ & 24153, 25314, 31524, 35142, 41352, 42513 \\
\end{tabular}
\end{center}

Simple permutations have recently had important applications in the
study of pattern closed classes of permutations \cite{AA02}.

Let $s_n$ denote the number of simple permutations of length $n$. We
shall be concerned with properties of the sequence $(s_n)$. Consider
the ordinary generating functions:
\begin{eqnarray*}
F(x) &=& \sum_{k=1}^{\infty} k! x^k ; \\ 
S(x) &=& \sum_{k=4}^{\infty} s_k x^k.
\end{eqnarray*}
We start $S(x)$ from $x^4$ because simple permutations of length $1$ and $2$ need 
special treatment.
Later in this section we will see that the coefficients of $S$ differ from those of
$-F^{\langle -1 \rangle}$ (functional inverse, not reciprocal) alternately by $2$ and
$-2$. The coefficients of $F^{\langle -1 \rangle}(x)$ were considered by Comtet
\cite[p. 171]{Co01} without any combinatorial
interpretation. The sequence of absolute values of these coefficients
appears as sequence A059372 of \cite{S01}, and the first few terms are:
\[
1,\, 2,\, 2,\, 4,\, 4,\, 48,\, 336,\, 2928,\, 28144,\, 298528,\, 3454432,\, 43286528.
\]
So we shall see that the numbers $s_n$ are:
\[
1,\, 2,\, 0,\, 2,\, 6,\, 46,\, 338,\, 2926,\, 28146,\, 298526,\, 3454434,\, 43286526.
\]
In
section \ref{nonPRecursiveness} we shall prove that $(s_n)$ is not
P-recursive (it cannot be defined by a linear recurrence with polynomial
coefficients).  In section \ref{asymptotics} we derive the asymptotic
behaviour of $s_n$ (the main term is $n!/\e^2$) and section
\ref{congruences} gives various congruences satisfied by the numbers
$s_n$.

In the remainder of this section we derive a structure theorem that
shows how arbitrary permutations are built from simple ones, and read
off from it equations satisfied by generating functions.  We begin
with some terminology and notation that will be used throughout.

A {\em block decomposition} of a permutation $\sigma$ is a partition
of $\sigma$ into blocks.  Of course, if $\sigma$ is simple there will
only be the two trivial block decompositions.  An example of a
non-trivial decomposition is $\sigma=67183524$ with blocks
$(67)(1)(8)(3524)$.

\begin{figure}
\begin{center}
\includegraphics[width=3cm]{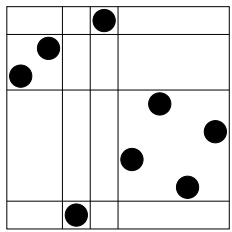}
\end{center}
\caption{A block decomposition of 67183524. The pattern of the block
decomposition is the permutation whose graph is defined by the
occupied cells, namely $3142$. Within each occupied cell, the
individual blocks also define permutations namely 12, 1, 1, and 2413.}
\label{BlockDecomposition}
\end{figure}

Given a block decomposition of $\sigma$, its \emph{pattern} is the
permutation defined by the relative order of the blocks. In the
example above, the pattern of the block decomposition
$(67)(1)(8)(3524)$ is $3142$. We may think of the permutation
$67183524$ as being constructed from the permutation $3142$ by {\em
inflating} each of the elements into a block, in this case the blocks
$12$, $1$, $1$, and $2413$ (we view each block as a permutation in its
own right). We write:
\[
67183524 = (3142)[12, 1, 1, 2413].
\]
This example is further illustrated in Figure
\ref{BlockDecomposition}.  The inflation procedure is an instance of
the wreath product for permutations \cite{AS01}.

A permutation which cannot be written in the form $(12)[\alpha,
\beta]$ is called {\em plus indecomposable}, and one which cannot be
written in the form $(21)[\alpha, \beta]$ is called {\em minus
indecomposable}.  Let $i_n$ denote the number of plus indecomposable
permutations of length $n$. The number of minus indecomposable
permutations of length $n$ is also $i_n$ as is easily seen by
considering the bijection on permutations of length $n$ which sends
$\pi$ to $\pi'$ where $\pi'(t) = n + 1 - \pi(t)$.

\begin{theorem}
\label{StructureTheorem}
For every non-singleton permutation $\pi$ there exists a unique simple 
non-singleton permutation $\sigma$, 
and permutations $\seq{\alpha}{k}$ such that
\[
\pi = \sigma[\seq{\alpha}{k}].
\]
Moreover, if $\sigma \neq 12, 21$ then $\seq{\alpha}{k}$ are also
uniquely determined. If $\sigma = 12$ (respectively $21$) then
$\alpha_1$ and $\alpha_2$ are uniquely determined subject to the
additional condition that $\alpha_1$ be plus (respectively minus)
indecomposable.
\end{theorem}

The caveat added for the case where $\sigma = 12$ (or 21) is
necessary, as is easily seen by considering $\pi = 123$. This can be
decomposed as $(12)[1,12]$ or as $(12)[12, 1]$. However, only the
former decomposition has a plus indecomposable first part.

\begin{proof}
We first of all suppose that $\pi$ has two distinct maximal proper
blocks $A$ and $B$ that have a non-empty intersection.  Then, as the
union of intersecting segments is a segment and the union of
intersecting ranges is a range, $A\union B$ is a block.  Because of
the maximality, $A\union B=[n]$.  But it is also clear that $A$ cannot
be an interior segment of $[n]$ nor can it define an interior range.
In other words we have
\[
\pi = \sigma[\alpha, \beta]
\]
where $\sigma = 12$ or $\sigma = 21$. These two possibilities are
obviously mutually exclusive.  In either case consider all
decompositions of $\pi$ as $\sigma[\gamma,\delta]$.  The intersection
of their $\gamma$ parts is also the $\gamma$ part of a decomposition
of this type.  So there is a unique such decomposition with smallest
$\gamma$ part.  Clearly, this part is plus indecomposable in the case
$\sigma=12$ and minus indecomposable if $\sigma=21$.

We next suppose that every pair of distinct maximal blocks has empty
intersection.  Obviously, then the maximal blocks form a block
decomposition of $\pi$ and this decomposition must be coarser than
every other block decomposition of $\pi$.  It follows that this
decomposition is the only one whose pattern $\sigma$ is simple and so
we obtain the unique representation claimed for $\pi$.
\end{proof}

We shall shortly see that this theorem gives relations between the
following three generating functions:
\begin{eqnarray*}
F(x) &=& \sum_{k=1}^{\infty} k! x^k ; \\ 
I(x) &=& \sum_{k=1}^{\infty} i_k x^k ; \\ 
S(x) &=& \sum_{k=4}^{\infty} s_k x^k.
\end{eqnarray*}
Note that our generating functions are all taken to have zero constant
term. This slightly unconventional choice turns out to be
algebraically convenient at several points.

From Theorem \ref{StructureTheorem} it is easy to see that there is a
one to one correspondence between the collection of all permutations
with length at least 2
and the collection of sequences:
\[
( \sigma, \, \seq{\alpha}{k} ).
\]
Here $\sigma$ may be any simple permutation of length $k\ge 2$, and if
$\sigma \neq 12, 21$ then $\alpha_1$ through $\alpha_k$ are arbitrary
permutations, while if $\sigma = 12$ (respectively 21), $\alpha_1$ is
plus-indecomposable (respectively minus indecomposable) and $\alpha_2$
is arbitrary.

This correspondence, together with the earlier observation that the
numbers of plus and minus indecomposable permutations of length $n$
are the same, translates naturally into the following equation:
\begin{equation}
\label{MainEquation}
F(x) = x + 2 I(x) F(x) + (S \circ F)(x).
\end{equation}
However, since a plus indecomposable permutation cannot correspond to
a sequence beginning with $12$, while all other sequences do represent
plus indecomposables, it is also clear from the correspondence that
\[
I(x) = x + I(x) F(x) + (S \circ F)(x).
\]
Solving this latter equation for $I$, and then substituting in
equation (\ref{MainEquation}) before solving for $S \circ F$ gives:
\[
(S \circ F)(x) = \frac{F(x) -F(x)^2}{1 + F(x)} - x.
\]
Now letting $t = F(x)$ we obtain:
\begin{equation}
\label{SEquation}
S(t) = t - \frac{2 t^2}{1 + t} - F^{\langle -1 \rangle}(t).
\end{equation}
We can also obtain an equation for the ordinary generating function of
plus indecomposable permutations through the observation that every
permutation decomposes into a sequence of plus indecomposable
permutations so
\[
F(x) = \frac{I(x)}{1 - I(x)}
\]
or equivalently
\begin{equation}
\label{IEquation}
I(x) = \frac{F(x)}{1 + F(x)}.
\end{equation}

Denoting the coefficient of $t^n$ in $F^{\langle -1 \rangle}(t)$ by
$\Com_n$ (in reference to Comtet who initiated the consideration of
this sequence) we obtain directly from equation (\ref{SEquation}) the
simple relationship that for $n \geq 4$:
\[
s_n = -\Com_n + (-1)^{n+1} \cdot 2.
\]

\section{Non P-recursiveness}\label{nonPRecursiveness}

A sequence of numbers $(a_n)$ is called P-recursive if it satisfies a
linear recurrence with polynomial coefficients. A power series is
called D-finite if it satisfies a linear differential equation with
polynomial coefficients. A sequence $(a_n)$ is P-recursive if and only
if its ordinary generating function $A(x)=\sum_{n}a_nx^n$ is
D-finite. More information on D-finiteness and P-recursiveness can be
found in Stanley \cite[Chapter 6]{St02}. If $a_n=n!$ then
$a_n-na_{n-1}=0$, and thus the sequence $(n!)$ is P-recursive. We show
that on the other hand neither sequence $(i_n)$ nor $(s_n)$ is
P-recursive.  By (\ref{SEquation}), instead of the latter sequence we
can work with $(\Com_n)$.

\begin{proposition}
\label{ISDiffEq}
The power series $I(x)$ and $C(x) = F^{\langle -1 \rangle}(x) = \sum_{k=1}^{\infty}
\Com_k x^k$ satisfy the differential equations
\begin{eqnarray*}
I' &=& -x^{-2}I^2+(x^{-2}+x^{-1})I-x^{-1} ; \\ C' &=&
\frac{C^2}{x-(1+x)C}.
\end{eqnarray*}
\end{proposition}

\begin{proof}
It follows from the recurrence for $n!$ that $F(x)$ satisfies
$x+xF+x^2F'=F$. Thus $F'=((1-x)F-x)/x^2$. Combining this with $F
=I/(1-I)$ we obtain the differential equation for $I(x)$. Similarly,
$C'=1/F'(C)= C^2/((1-C)x-C)$ which is the differential equation for
$C(x)$.
\end{proof}

Klazar \cite{Kl01} used the following method to show that a
sequence $(a_n)$ is not P-recursive. Suppose that the ordinary
generating function $A(x)$ is non-analytic and satisfies a first order
differential equation $A'=R(x,A)$ where $R$ is some
expression. Differentiating this relationship and replacing $A'$ by
$R(x,A)$, the derivatives of $A$ are expressed as $A^{(k)}=R_k(x,A)$;
$R_0(x,A)=A$ and $R_1(x,A)=R(x,A)$. Substituting $R_k(x,A)$ in the
equation of D-finiteness
\[
b_0A+b_1A'+b_2A''+\cdots+b_sA^{(s)}=0,
\]
where $s\ge 1$, $b_i\in\C(x)$ and $b_s\ne 0$, we get a
non-differential equation $\sum_{k=0}^s b_kR_k(x,A)=0$.  If $R$ is
such that the expressions $R_0,R_1,R_2,\dots$ are (i) analytic or even
algebraic and (ii) linearly independent over $\C(x)$, we have a
nontrivial analytic equation for $A$. This implies that $A$ is
analytic (see Klazar's paper \cite{Kl01} for more details) which is a
contradiction. So $A$ cannot be D-finite and the sequence of its
coefficients cannot be P-recursive.

To state the result of \cite{Kl01} precisely, we remind the reader
that a power series $R(x,y)\in\C[[x,y]]$ is analytic if it absolutely
converges in a neighborhood of the origin and that
$R(x,y)\in\C((x,y))$ is an analytic Laurent series if, for some
positive integer $k$, $(xy)^kR(x,y)\in\C[[x,y]]$ is analytic. Theorem
1 of \cite{Kl01} says that if $A\in\C[[x]]$ is non-analytic,
$R(x,y)\in\C((x,y))$ is analytic, $A'=R(x,A)$, and $R$ contains at
least one monomial $ax^iy^j$, $a\ne 0$, with $j<0$, then $A$ is not
D-finite. This result applies directly neither to $I(x)$ nor $C(x)$
(see Proposition~\ref{ISDiffEq}) because in the case of $I(x)$ the
last condition on $R$ is not satisfied and in the case of $C(x)$ the
right hand side $R$ even cannot be expanded as a Laurent series.

However, the substitution $x-(1+x)C(x)=\theta(x)$ transforms the
second differential equation of Proposition~\ref{ISDiffEq} into
\[
\theta'=-\frac{x^2}{1+x}\cdot\frac{1}{\theta}+\frac{1+2x}{1+x}.
\]
Now all conditions are satisfied ($F(x)$ is clearly non-analytic which
implies that $C(x)$ and $\theta(x)$ are non-analytic) and thus
$\theta(x)$ is not D-finite by Theorem 1 of \cite{Kl01}. The
dependence of $C(x)$ and $S(x)$ on $\theta(x)$ and the fact that
D-finite power series form a $\C(x)$-algebra (\cite[Theorem
6.4.9]{St02}) shows that neither $C(x)$ nor $S(x)$ is D-finite.

In order to deal with the case of $I(x)$, we use this opportunity to
complement Theorem 1 of \cite{Kl01} in which $R\in\C((x,y))$ by the
following theorem which treats the case $R\in\C(x,y)$. Neither of the
theorems subsumes the other because not every rational function in $x$
and $y$ can be represented by an element of $\C((x,y))$ (as we have
seen) and, of course, not every Laurent series sums up to a rational
function. However, the next theorem seems to be more useful because in
both examples in \cite{Kl01} and both examples here the right hand
side $R(x,y)$ is, in fact, a rational function.

\begin{theorem}
\label{nonPrec}
Let $P,Q\in\C[x,y]$ be two nonzero coprime polynomials and
$A\in\C[[x]]$ be a non-analytic power series which satisfies the
differential equation
$$
A'=\frac{P(x,A)}{Q(x,A)}.
$$
If $\deg_y Q=0$ and $\deg_y P\le 1$ then $A$ is, trivially,
D-finite. In all remaining cases $A$ is not D-finite.
\end{theorem}

\begin{proof}
The first claim is clear. If $\deg_y Q=0$ and $r=\deg_y P\ge 2$ then
$A'=a_0+a_1A+\cdots+a_rA^r$ where $a_i\in\C(x)$, $r\ge 2$, and $a_r\ne
0$. Differentiation by $x$ gives
$$
A^{(k)}=R_k(x,A)=a_{0,k}+a_{1,k}A+\cdots+a_{kr-k+1,k}A^{kr-k+1}
$$
where $a_{i,j}\in\C(x)$ and
$$
a_{kr-k+1,k}=r(2r-1)(3r-2)\dots((k-1)r-k+2)a_r^k\ne 0.
$$
Thus $R_k(x,y)\in\C(x)[y]$ have $y$-degrees $kr-k+1$, $k=0,1,2,\dots$,
which is for $r\ge 2$ a strictly increasing sequence. Therefore
$R_0,R_1,R_2,\dots$ are linearly independent over $\C(x)$ and, by the
above discussion, $A$ is not $D$-finite.

In the remaining case $\deg_y Q\ge 1$. Differentiation of
$A'=R(x,A)=P(x,A)/Q(x,A)$ by $x$ gives $A^{(k)}=R_k(x,A)$ where
$R_k(x,y)\in\C(x,y)$. For example,
\begin{eqnarray*}
R_2&=&\frac{(P_x+P_yR_1)Q-P(Q_x+Q_yR_1)}{Q^2}\\
&=&\frac{P_xQ-PQ_x}{Q^2}+\frac{P(P_yQ-PQ_y)}{Q^3}.
\end{eqnarray*}
Let $\alpha$, $Q(x,\alpha)=0$, be a pole of $R_1(x,y)$ of order
$\ord_{\alpha}(R_1)=\ord_{\alpha}(P/Q)=-\ord_{\alpha}(Q)=l\ge 1$. We
have $\ord_{\alpha}((P_xQ-PQ_x)Q^{-2})\le 2l$ and
$\ord_{\alpha}(P(P_yQ-PQ_y)Q^{-3})=3l+\ord_{\alpha}(P_yQ-PQ_y)=2l+1$
since $\ord_{\alpha}(P)=0$, $\ord_{\alpha}(P_yQ)\le -l$, and
$\ord_{\alpha}(PQ_y)=-l+1$. So $\ord_{\alpha}(R_2)=2l+1$.  In general,
the same argument shows that
$\ord_{\alpha}(R_{k+1})=2\cdot\ord_{\alpha}(R_k)+1$.  Hence
$\ord_{\alpha}(R_k)=2^{k-1}l+2^{k-1}-1$, $k=1,2,\dots$.  This is a
strictly increasing sequence and we conclude again, since
$R_0,R_1,R_2,\dots$ are linearly independent over $\C(x)$, that $A$ is
not $D$-finite.
\end{proof}

Proposition~\ref{ISDiffEq} and Theorem~\ref{nonPrec} show that $I(x)$
is not $D$-finite and we can summarize the results of this section in
the following corollary.

\begin{corollary}
The sequences $(i_n)$, $(\Com_n)$, and $(s_n)$ are not P-recursive.
\end{corollary}

\section{Asymptotics}\label{asymptotics}

We turn now to the computation of an asymptotic expansion for the
numbers $s_n$. We will prove that:
\begin{theorem}
\[
s_n = \frac{n!}{\e^2} \left( 1 - \frac{4}{n} + \frac{2}{n(n-1)} +
O(n^{-3}) \right) .
\]
\end{theorem}

Our methods are such that, in principle, higher order terms could be
obtained as a matter of brute force computation. In order to carry out
this expansion we will first consider permutations which may not be
simple, but whose non-trivial blocks all have length greater than some
fixed value $m$. We will apply inclusion-exclusion arguments
(dressed in the form of generating functions 
\cite{FS01, GJ01}), an argument which
allows us to reduce the number of terms considered, and a
bootstrapping approach. 

The case $m = 2$, was already considered by Kaplansky
\cite{Ka01}. Permutations of this type are those in which no two
elements consecutive in position are also consecutive in value (in
either order). These were called irreducible permutations by Atkinson
and Stitt
\cite{AS01}, but there is no standard terminology in the field.
Indeed the permutations that we have referred to as plus and minus
indecomposable have also been called irreducible in other contexts.

An amusing equivalent form for the case $m = 2$ is that the number of
such permutations is also the number of ways of placing $n$ mutually
non-attacking {\em krooks} on an $n \times n$ chessboard. A krook is a
piece which can move either like a king, or a rook in
chess. Kaplansky's expansion is:
\[
\frac{n!}{\e^2} \left( 1 - \frac{2}{n(n-1)} + O(n^{-3}) \right).
\]

In fact he derives asymptotic forms for the number of permutations
containing exactly $r$ blocks of length 2 for any $r$. Our methods
parallel his, and could also be used to derive such detailed
information.

The decomposition provided by Theorem \ref{StructureTheorem} of a
permutation into its maximal proper blocks represents a top down view
of how non-simple permutations are constructed from simple ones. There
is a corresponding bottom-up view that focuses on minimal blocks, put
together in an arbitrary order. By a {\em minimal block} in $\pi$ we
mean a non-singleton block in $\pi$ minimal with respect to
inclusion. Note that the pattern of each minimal block is that of a
simple permutation. Any permutation can be decomposed into minimal
blocks and singletons, e.g., $3524716=(3524)(7)(1)(6)$. However, this
decomposition is not unique, for two essentially different reasons.
The first one is that decompositions $\pi = \sigma[\seq{\alpha}{k}]$,
where $\sigma$ is arbitrary and $\alpha_i$ are simple, are not unique
because it may be possible to coalesce singletons into simple blocks,
or vice versa. Thus besides $3524716=2413[2413,1,1,1]$ we also have
$3524716=3524716[1,1,1,1,1,1,1]$. The second problem is that we
require any two minimal blocks to be disjoint. While this is
necessarily true whenever either of them has length more than 2, two
minimal blocks of length 2 may intersect, as in $123$. Thus we
consider decompositions $\pi = \sigma[\seq{\alpha}{k}]$ where $\sigma$
is arbitrary and each $\alpha_i$ is either 1, a simple permutation of
length at least 4, or the identity permutation of length at least 2 or
its reverse.  We refer to blocks of the latter type as {\em clusters} in $\pi$.

By using clusters we have solved the second problem but the
non-uniqueness remains and, moreover, we have introduced another
source of it: consecutive (reversed) identical permutations may
coalesce into longer (reversed) identical permutations, as in
$345612=21[1234,12]=231[12,12,12]$. To remedy the non-uniqueness we
introduce the notion of {\em marking} a permutation. A {\em marked
permutation} $(\pi, M)$ consists of a permutation $\pi$ and a
collection $M$ of minimal blocks of $\pi$. A {\em marked cluster} in
$(\pi, M)$ is a maximal chain of marked overlapping minimal blocks of
length 2 (a marked cluster may be a proper subset of a maximal cluster).  Let
${\cal B}_1$ denote the set of all simple permutations of length at
least 4 and ${\cal B}_2$ denote the set of all identical permutations
of length at least 2 and their reversals. Marking makes our
decomposition unique:

\begin{theorem}
\label{BottomUpStructureTheorem} 
Let $X$ be the set of all marked permutations $(\pi, M)$ and $Y$ be
the set of all sequences $(\sigma;\seq{\alpha}{k})$ where $\sigma$ is
any permutation of length $k\ge 1$ and $\alpha_i\in\{1\}\cup{\cal
B}_1\cup{\cal B}_2$.  There is a bijection between the sets $X$ and
$Y$ such that if $(\pi, M)\mapsto (\sigma;\seq{\alpha}{k})$, where $r$
of the $\alpha_i$ belong to ${\cal B}_1$ and $s$ of them to ${\cal
B}_2$, then
\[
\pi = \sigma[\seq{\alpha}{k}]
\]
and $|M|=r+l-s$ where $l$ is the total length of the $\alpha_i$ belonging to 
${\cal B}_2$.
\end{theorem}

\begin{proof}
Given a marked permutation, collapse its marked minimal blocks of
length at least 4 and its marked clusters into singletons. This gives
the permutation $\sigma$.  If the $i$-th term of $\sigma$ was not
obtained by collapse then $\alpha_i=1$, otherwise $\alpha_i$ equals to
the corresponding element of ${\cal B}_1\cup{\cal B}_2$. Since each
$\alpha_i\in{\cal B}_1$ contributes 1 to $|M|$ and each
$\alpha_i\in{\cal B}_2$ of length $m$ contributes $m-1$, we have
$|M|=r+l-s$. It is clear that $\pi = \sigma[\seq{\alpha}{k}]$ and that
$(\pi, M)$ can be uniquely recovered from $(\sigma;\seq{\alpha}{k})$.
\end{proof}

Now suppose $m$ to be some fixed value (we will later make choices of
$m$ suitable for our purposes, but will always assume that $m \geq 2$
since smaller values of $m$ are trivial). Each permutation $\pi$ has
an associated collection $B_m(\pi)$ consisting of the minimal
blocks of $\pi$ whose length is less than or equal to
$m$. So, if $\pi$ is simple and of length greater than $m$, $B_m(\pi)$
is empty, while for $\pi = 5672413$, $B_2(\pi) = \set{56, 67}$, and
$B_4(\pi) = \set{56, 67, 2413}$.  An {\em $m$-marking} of
$\pi$ is simply a subset of $B_m(\pi)$. We consider the generating
function:
\[
F_m(x, v) = \sum_{\pi} x^{|\pi|} \sum_{M \ci B_m(\pi)} v^{|M|} =
\sum_{\pi} x^{|\pi|} (1+v)^{|B_m(\pi)|}.
\]
Then of course $F_m(x, -1)$ is the ordinary generating function for
permutations all of whose non-singleton blocks have length greater than $m$.

We remark that $F_m(x, t-1)$ is the generating function where the
coefficient of $x^n t^k$ is precisely the number of permutations of
length $n$ with $k$ minimal blocks of length less than or equal to
$m$. 

Let
\[
S_m(x) = \sum_{j=4}^m s_j x^j.
\]
We apply the bijection of Theorem~\ref{BottomUpStructureTheorem} to
marked permutations which contain no marked minimal blocks of length
more than $m$. It follows that the generating function of the
corresponding permutations $\alpha\in\{1\}\cup{\cal B}_1\cup{\cal
B}_2$, in which $x$ counts the length and $v$ the contribution to
$|M|$, is

\[
x+vS_m(x)+\frac{2vx^2}{1-vx}.
\]
So:
\[
F_m(x, v) = \sum_{k \geq 1} k! \, \left( x + \frac{2 v x^2}{1 - v x} +
v S_m(x) \right)^k
\]
from which it follows that:
\begin{equation}
\label{MAsympEq}
f_m(x) := F_m(x,-1) = \sum_{k \geq 1} k! \, \left(x - \frac{2 x^2}{1 +
x} - S_m(x) \right)^k.
\end{equation}

Before using this equation to derive asymptotic information about
$s_n$ we digress briefly to show how it can be used to obtain an
alternative derivation of (\ref{SEquation}). Instead of using
$S_m(x)$ in (\ref{MAsympEq}), use $S(x)$. This gives us $f_{\infty}(x)$, an
ordinary generating function for permutations having no minimal
block. The only such permutation is $1$ so
$f_{\infty}(x) = x$. That is:
\[
x = F(x - \frac{2 x^2}{1 +
x} - S(x))
\]
which yields (\ref{SEquation}) after applying $F^{\langle -1 \rangle}$
to both sides.

Now recall that $f_m(x)$ is the generating function for permutations all
of whose blocks have length greater than $m$. In order to make use of
these generating functions in the asymptotic computation of $s_n$ we
must determine a suitable value of $m$ so that $f_m$ provides useful
information about $s_n$. To that end the following lemma is useful.

\begin{lemma}
If $p_{n,k}$ denotes the number of permutations of length $n$ which
contain a minimal block of length $k$ then for any fixed positive
integer $c$:
\[
\sum_{k = c+2}^{n-c} \frac{p_{n,k}}{n!}  = O(n^{-c}).
\]
\end{lemma}

\begin{proof}
First observe that
\[
p_{n,k} \leq s_k (n-k+1) (n-k+1)!
\]
since the right hand side counts the number of ways to choose the
structure of a block of length $k$, to choose its minimal element, and
to arrange it with other elements, so it overcounts permutations with
more than one such block.

The estimate given then follows directly by using the fact that $s_k
\leq k!$. Only the two extreme terms in the sum can have magnitude as
large as $O(n^{-c})$, and the remaining terms have magnitude
$O(n^{-c-1})$. Since there are fewer than $n$ terms, the result
follows.
\end{proof}

So when seeking an asymptotic expansion of $s_n/n!$ with an error term
of $O(n^{-c-1})$ we may count instead the permutations which contain
no blocks of length less than or equal to $c+2$, or greater than or
equal to $n-c$. In particular, as a direct consequence of the result
quoted above due to Kaplansky \cite{Ka01} we obtain:

\begin{observation}
\label{FirstOrder}
\[
\frac{s_n}{n!} = \frac{1}{\e^2} + O(n^{-1}).
\]
\end{observation}

An alternative proof of this result follows from a more general
theorem of Bender and Richmond \cite{BR01} which provides the first
order asymptotics of a class of series which include the inverse
series of $F(x)$.

We will set as our goal to obtain the asymptotics of $s_n/n!$ with
error term $O(n^{-3})$. However, the technique we use is completely
general, and could be applied, at the expense of a great deal of
tedious computation, to any fixed error bound of this type. By the
remarks above, we may ignore minimal block sizes between $5$ and $n-3$
inclusive. We first consider $f_4(x)$ which enumerates permutations
having no minimal blocks of size less than or equal to 4. Recall that:
\[
f_4(x) = \sum_{k \geq 1} k! \, \left(x - \frac{2 x^2}{1 + x} - 2 x^4
\right)^k.
\]
So:
\begin{eqnarray}
\frac{1}{n!} [t^n] f_4(t) &=& \frac{1}{n!} \sum_{k=0}^{\infty} k! \,
[t^n] \left( t - \frac{2t^2}{1 + t} - 2t^4 \right)^k \nonumber \\ 
{} &=& \frac{1}{n!}  \sum_{k=0}^{\infty} k! \, [t^{n-k}] \left( 1 -
\frac{2t}{1 + t} -2 t^3 \right)^k \nonumber \\ 
{} &=& \frac{1}{n!} \sum_{l = 0}^n (n-l)! \, [t^l] \left( 1 -
\frac{2t}{1 + t} - 2 t^3 \right)^{n-l} \nonumber \\ 
{} &=& \frac{1}{n!} \sum_{l = 0}^n (n-l)! \sum_{i = 0}^{l} (-2)^i {n-l
\choose i} [t^l] \left( \frac{t}{1 + t} + t^3 \right)^{i} \nonumber \\ 
{} &=& \frac{1}{n!} \sum_{l = 0}^n (n-l)! \sum_{i = 0}^{l} (-2)^i {n-l
\choose i} [t^{l-i}] \left( \frac{1}{1 + t} + t^2 \right)^{i}.
\label{F4Equation}
\end{eqnarray}

Consider now any fixed value of $l$ in equation (\ref{F4Equation}). In
order to obtain terms whose order in $n$ is $n^{-2}$ or more, we need
only consider the values $l-2\leq i\leq l$. Despite the fact that we
sum over values of $l$ running from $0$ through $n$, we may safely
ignore the other terms. As we shall see in computing the three
significant terms the summation over $l$ does not affect the order of
the terms.

So, the three terms that we need to consider are:
\begin{equation}
\begin{array}{l}
\frac{(n-l)!}{n!}  (-2)^l {n-l \choose l} + \\ 
\frac{(n-l)!}{n!} \left( (-2)^{l-1} {n-l \choose l-1} (-l+1) \right) + \\
\frac{(n-l)!}{n!} \left( (-2)^{l-2} {n-l \choose l-2} \left( (-l +
2)(-l + 1)/2 + l-2 \right) \right).
\end{array}
\end{equation}
Each of these terms will be converted to the form:
\[
\frac{(-2)^l}{l!} \left( \mbox{an asymptotic expansion in $n$}
\right).
\]
Since the first two and the first part of the third, are the same as
those arising in the $m = 2$ case, we can make use of their known
form, that is, use the asymptotics from Kaplansky's result, leaving
only the term
\begin{eqnarray*}
\frac{(n-l)!  (-2)^{l-2} (l-2)}{n!} {n-l \choose l-2} &=&
\frac{(-2)^l}{l!} \left( \frac{l(l-1)(l-2)}{4} \, \frac{(n-l)!
(n-l)!}{n!  (n-2l+2)!} \right) \\ {} &=& \frac{(-2)^l}{l!} \left(
\frac{l(l-1)(l-2)}{4n(n-1)} + O(n^{-3}) \right)
\end{eqnarray*}
Summing this expression over $l$ gives $-2 \e^{-2}/n(n-1) + O(n^{-3})$.

Now we combine this additional term with Kaplansky's results to give
the asymptotic expansion of $[t^n] f_4(t)$ through three terms as:
\[
[t^n] f_4(t) = \frac{n!}{\e^2} \left( 1 - \frac{4}{n(n-1)} + O(n^{-3})
\right).
\]

Finally we use this in establishing the second order asymptotics of
$s_n$. From Observation \ref{FirstOrder} applied to $s_{n-1}$ we
obtain:
\[
s_{n-1} = \frac{n!}{\e^2} \left( \frac{1}{n} + O(n^{-2}) \right).
\]
Furthermore, the number of permutations of length $n$ containing a
simple block of length $n-1$ is precisely $4 s_{n-1}$. Since, in
computing the $1/n$ term in the expansion of $s_n$ we can ignore
contributions arising from blocks of length $n-2$, and since the
events of having a simple block of length from 2 to 4, and having a
simple block of length $(n-1)$ are disjoint:
\begin{eqnarray*}
s_n &=& [t^n] f_4(t) - 4 s_{n-1} + O(n^{-2} n!) ; \\ {} &=&
\frac{n!}{\e^2} \left( 1 - \frac{4}{n} + O(n^{-2}) \right).
\end{eqnarray*}

We apply this bootstrap approach once more to get the second order
behaviour. We now know that:
\begin{eqnarray*}
s_{n-1} &=& \frac{n!}{\e^2} \left( \frac{1}{n} - \frac{4}{n(n-1)} +
O(n^{-3}) \right) \\ s_{n-2} &=& \frac{n!}{\e^2} \left(
\frac{1}{n(n-1)} + O(n^{-3}) \right).
\end{eqnarray*}
Furthermore there are $18 s_{n-2}$ permutations of length $n$
containing a simple block of length $n-2$. However, of these $8
s_{n-2}$ also contain a simple block of length $2$. So:
\begin{eqnarray*}
s_n &=& [t^n] f_4(t) - 4 s_{n-1} - 10 s_{n-2} + O(n^{-3} n!) \\ {} &=& \frac{n!}{\e^2}
 \left( 1 - \frac{4}{n} + \frac{2}{n(n-1)} + O(n^{-3}) \right),
\end{eqnarray*}
as we claimed at the beginning of this section.

Finally, in this section we note that the asymptotic estimate of $s_n$
is, as might be expected, a poor approximation.  For example,
$s_{20}=264111424634864638$ and our asymptotic estimate has a relative
error of about $3.89\times 10^{-3}$.

\section{Congruences}\label{congruences}

In this section we derive congruence properties of the numbers
$\mathrm{Com}_n$ for the moduli $2^a$ and $3$ (from which follow similar
congruences for $s_n$).  Our main tool is the following result that
follows immediately from the Lagrange inversion formula.

\begin{lemma}\label{lagrange}
$$
n\cdot\mathrm{Com}_n=[x^{n-1}]\left(\sum_{k\ge
0}(-1)^k(2!x+3!x^2+\cdots)^k\right)^n.
$$
\end{lemma}

For a prime $p$, let $\ord_p(n)$ denote the largest integer $m$ such
that $p^m$ divides $n$. As the following table shows,
$\ord_2(\mathrm{Com}_n)$ is unexpectedly large:

\bigskip\begin{center}
\begin{tabular}{l||r|r|r|r|r|r|r|r|r|r|r|r|r|r|r}
$n$ & $1$ & $2$ & $3$ & $4$ & $5$ & $6$ & $7$ & $8$ & $9$ & $10$ &
$11$ & $12$ & $13$ & $14$ & $15$\\ \hline$\ord_2(\mathrm{Com}_n)$ &
$0$ & $1$ & $1$ & $2$ & $2$ & $4$ & $4$ & $4$ & $4$ & $5$ & $5$ & $15$
& $13$ & $12$ & $12$
\end{tabular}

\bigskip
\begin{tabular}{r|r|r|r|r|r|r|r|r|r|r|r|r|r|r}
$16$ & $17$ & $18$ & $19$ & $20$ & $21$ & $22$ & $23$ & $24$ & $25$ &
$26$ & $27$ & $28$ & $29$ & $30$ \\ \hline $8$ & $8$ & $9$ & $9$ &
$10$ & $10$ & $12$ & $12$ & $14$ & $14$ & $15$ & $15$ & $17$ & $17$ &
$22$
\end{tabular}

\end{center}

\bigskip\noindent In Theorem~\ref{dolniodhad} we give a lower bound on
$\ord_2(\mathrm{Com}_n)$ which is tight for infinitely many $n$ and we
completely characterize the values of $n$ for which the equality is
attained.

For convenience we note the following result that follows directly
from the well-known formula
\[
\ord_p(m!)=\left\lfloor\frac{m}{p}\right\rfloor+\left\lfloor\frac{m}{p^2}\right\rfloor+\cdots
\]

\begin{lemma}\label{faktorial}
For all $m$, $\ord_2((m+1)!)\ge\left\lceil\frac{m}{2}\right\rceil$
where equality holds if and only if $m=1$ or $2$. Also, $\ord_3(m!)\le
m-1$ for all $m$.
\end{lemma}

\begin{theorem}\label{dolniodhad}
Let $m=\lfloor n/2\rfloor$. Then
$$
\ord_2(\mathrm{Com}_n)\ge\left\lceil\frac{n-1}{2}\right\rceil.
$$
Equality holds if and only if ${3m\choose m}$ is odd and this happens
if and only if the binary expansion of $m$ has no two consecutive unit
digits.
\end{theorem}

\begin{proof}
Let the numbers $b_k$, $k\ge 0$, be defined by
$$
\sum_{k\ge 0}b_kx^k=\sum_{k\ge 0}(-1)^k(2!x+3!x^2+\cdots)^k.
$$
Thus $b_0=1$ and for $k\geq 1$,
$$
b_k= \sum_{\stackrel{\scriptstyle c_1,c_2,\dots,c_s\ge
1}{c_1+c_2+\cdots+c_s=k}}(-1)^s\cdot(c_1+1)!\cdot
(c_2+1)!\cdot\dots\cdot(c_s+1)!.
$$
By Lemma \ref{lagrange},
$$
n\cdot\mathrm{Com}_n= \sum_{\stackrel{\scriptstyle
k_1,k_2,\dots,k_n\ge 0}{k_1+k_2+\cdots+k_n=n-1}}b_{k_1}b_{k_2}\dots
b_{k_n}.
$$
By Lemma~\ref{faktorial}, $\ord_2((c+1)!)\ge c/2$ for all $c$. Hence,
for all $k$ and $n$,
$$
\ord_2(b_k)\ge\frac{k}{2}\ \mbox{ and }\
\ord_2(n\cdot\mathrm{Com}_n)\ge\frac{n-1}{2}.
$$
In particular, for odd $n$ we have
$\ord_2(\mathrm{Com}_n)=\ord_2(n\cdot\mathrm{Com}_n)\ge(n-1)/2$.

To obtain the more exact result of the theorem we need the following
better estimates for $\ord_2(b_k)$:$$ \ord_2(b_k)\left\{
\begin{array}{ll}
=k/2&\mbox{for even $k$;}\\ =(k+1)/2&\mbox{for $k\equiv 1$
mod $4$;}\\ >(k+1)/2&\mbox{for $k\equiv 3$ mod $4$.}
\end{array}
\right.
$$
To prove them we look more closely at the sum for $b_k$. Suppose first
that $k$ is even. Then the sum has exactly one summand with $\ord_2$
equal to $k/2$, namely that with $c_1=c_2=\dots=c_{k/2}=2$ (by
Lemma~\ref{faktorial}, $\ord_2((c+1)!)=c/2$ only if $c=2$), and the
other summands have $\ord_2$ bigger than $k/2$. Hence
$\ord_2(b_k)=k/2$. Now suppose that $k$ is odd. Then each summand has
an odd number of odd $c_i$'s. The summands $t$ with three and more odd
$c_i$'s satisfy $\ord_2(t)\ge(k+3)/2$ (each odd $c_i$ contributes
$1/2$ to $k/2$). The same is true if $t$ has only one odd $c_i$ but
that $c_i$ is not $1$ (by Lemma~\ref{faktorial},
$\ord_2((c+1)!)\ge(c+3)/2$ for odd $c>1$), or if some even $c_i$ is
not $2$ (Lemma~\ref{faktorial}).  The remaining summands $t$, in which
$c_i=2$ with multiplicity $(k-1)/2$ and once $c_i=1$, satisfy
$\ord_2(t)=(k+1)/2$. We see that, for odd $k$, $\ord_2(b_k)=(k+1)/2$
if and only if the number of the remaining summands is odd. This
number equals $(k-1)/2+1=(k+1)/2$. So $\ord_2(b_k)=(k+1)/2$ if and
only if $k\equiv 1$ mod $4$.

Let $n=2m+1$ be odd. If $s$ is a summand of the above sum for
$n\cdot\mathrm{Com}_n$, then $\ord_2(s)=(n-1)/2$ if and only if all
$k_i$ in $s$ are even; other summands $t$ have $\ord_2(t)>(n-1)/2$. It
follows that $\ord_2(\mathrm{Com}_n)=(n-1)/2$ if and only if the
number of the former summands $s$ is odd.  This number equals
$$
[x^{n-1}]\left(\sum_{r\ge
0}x^{2r}\right)^n=[x^{n-1}]\frac{1}{(1-x^2)^n}= [x^{n-1}]\sum_{r\ge
0}{n+r-1\choose r}x^{2r}={3m\choose m}.
$$

Let $n=2m$ be even. We know that $\ord_2(b_k)=k/2$ for even $k$ and
$\ord_2(b_k)\ge (k+1)/2$ for odd $k$. In the sum for
$n\cdot\mathrm{Com}_n$, every composition $k_1+k_2+\cdots+k_n=n-1$ of
$n-1$ has an odd number of odd parts. For any $t$-tuple
$l_1,l_2,\dots,l_t$, where $t$ and all $l_i$ are odd and
$l_1+\cdots+l_t\le n-1$, we let $S(l_1,l_2,\dots,l_t)$ denote the sum
of those $b_{k_1}b_{k_2}\dots b_{k_n}$ with $k_1+k_2+\cdots+k_n=n-1$
in which $k_i=l_i$, $1\le i\le t$, and $k_i$ is even for $i>t$. It
follows that
$$
n\cdot\mathrm{Com}_n=\sum {n\choose t}S(l_1,l_2,\dots,l_t)
$$
where we sum over all mentioned $t$-tuples $l_1,l_2,\dots,l_t$. By the
properties of $\ord_2$ and of the numbers $b_k$,
$\ord_2(S(l_1,l_2,\dots,l_t))\ge(n+t-1)/2$. Also, for odd $t$ we have
$\ord_2({n\choose t})=\ord_2(\frac{n}{t}{n-1\choose
t-1})=\ord_2(n)-\ord_2(t)+\ord_2({n-1\choose t-1})\ge\ord_2(n)$, and
$\ord_2({n\choose 1})=\ord_2(n)$. It follows that
$\ord_2(\mathrm{Com}_n)\ge n/2$ and, moreover,
$\ord_2(\mathrm{Com}_n)=n/2$ if and only if
$$
\ord_2\left(\sum_{l\le n,\;l\;\mathrm{odd}}S(l)\right)=n/2.
$$
In the last sum still many summands have $\ord_2$ bigger than $n/2$:
if $l\equiv 3$ mod $4$ then $\ord(S(l))>n/2$.  On the other hand, if
$l\equiv 1$ mod $4$ then each summand $b_lb_{k_2}\dots b_{k_n}$ in
$S(l)$ has $\ord_2(b_lb_{k_2}\dots b_{k_n})=n/2$. We conclude that
$\ord_2(\mathrm{Com}_n)=n/2$ if and only if the number $c(n)$ of
compositions of $n-1$ into $n$ parts, where the first part is $\equiv
1$ mod $4$ and the remaining $n-1$ parts are even (zero parts are
allowed), is odd. We have
\begin{eqnarray*}  
c(n)&=&[x^{n-1}]\frac{x}{1-x^4}\cdot\frac{1}{(1-x^2)^{n-1}}
=[x^{n-1}]\frac{x}{1+x^2}\cdot\frac{1}{(1-x^2)^n}\\
&\equiv&[x^{n-1}]\frac{x}{1-x^2}\cdot\frac{1}{(1-x^2)^n}=[x^{n-1}]\frac{x}{(1-x^2)^{n+1}}\
\mathrm{mod}\ 2\\ &=&{3m-1\choose m-1}\equiv \frac{3m}{m}{3m-1\choose
m-1}\ \mathrm{mod}\ 2\\ &=&{3m\choose m}.
\end{eqnarray*}
It was noted by Kummer \cite{kumm}, see also Singmaster
\cite{Si01}, that $\ord_p({a+b\choose b})$ is equal to the number of
carries required when adding $a$ and $b$ in the $p$-ary
notation. Applying this for $p=2$, $a=m$, and $b=2m$, we get the
stated criterion.
\end{proof}

\begin{corollary}\label{moc2}
For all $n\ge 3$,
$$
s_n\equiv\left\{
\begin{array}{rll}
2 & \mathrm{mod}\ 2^{(n-1)/2} & \mbox{ for odd $n$} ; \\ -2 &
\mathrm{mod}\ 2^{n/2} & \mbox{ for even $n$.}
\end{array}
\right.
$$ 
\end{corollary}

Let
$$
C_n=\frac{1}{n+1}{2n\choose n}
$$
be the $n$th Catalan number.

\begin{proposition}
For all $n$, $\mathrm{Com}_n\equiv C_{n-1}\ \mathrm{mod}\ 3$.
\end{proposition}

\begin{proof}
We have, for every non-negative integer $k$,
$$
(2!x+3!x^2+\cdots)^k=(2x)^k+3a_k(x)
$$
with $a_k(x)\in\Z[[x]]$.  Thus
\begin{eqnarray*}
\sum_{k\ge 0}(-1)^k(2!x+3!x^2+\cdots)^k&=& \frac{1}{1+2x}+3\sum_{k\ge
0}(-1)^ka_k(x)\\ &=&\frac{1}{1+2x}+3b(x)
\end{eqnarray*}
with $b(x)\in\Z[[x]]$. Let $m=\ord_3(n)$. Since $\ord_3(k!)\le k-1$
for every $k$ (Lemma~\ref{faktorial}), we have
$${\textstyle \ord_3\Big(3^k{n\choose k}\Big)\ge m+1\ \mbox{for
$k=1,2,\dots,n$.}  }
$$
By Lemma \ref{lagrange},
\begin{eqnarray*}
n\cdot\mathrm{Com}_n=[x^{n-1}]\left(\frac{1}{1+2x}+3b(x)\right)^n&\equiv&
[x^{n-1}]\frac{1}{(1+2x)^n}\ \mathrm{mod}\ 3^{m+1}\\
&=&(-2)^{n-1}{2n-2\choose n-1}.
\end{eqnarray*}
Canceling in the last congruence the common factor $3^m$, we get
$$
\frac{n}{3^m}\cdot\mathrm{Com}_n\equiv\frac{(-2)^{n-1}}{3^m}{2n-2\choose
n-1} \equiv\frac{1}{3^m}{2n-2\choose n-1}\ \mathrm{mod}\ 3.
$$
Since $n/3^m\not\equiv 0\ \mathrm{mod}\ 3$, we can divide by it and
get
$$
\mathrm{Com}_n\equiv\frac{1}{n}{2n-2\choose n-1}\ \mathrm{mod}\ 3.
$$
\end{proof}

\begin{corollary}\label{catalan}
For all $n>2$,
$$
s_n \equiv -C_{n-1}+(-1)^n\ \mathrm{mod}\ 3.
$$
\end{corollary}

\section{Concluding remarks}

The simplicity property for permutations does not seem to have been
studied until very recently \cite{Mu01, AA02}.  We have begun the study
of the numbers $s_n$ by showing that they are not P-recursive, giving
the first few terms of their asymptotic expansion, and showing that
they satisfy some unexpected congruence properties. 

These results suggest a number of natural continuations. Although, in
principle, we could obtain more terms of the asymptotic expansion the
entire expansion remains elusive, and computing it seems to be rather a
difficult problem. On the other hand we have some computational
evidence to suggest that the sequence $\Com_n$ has additional
congruence properties, particularly with respect to odd primes.

We suggest also some algorithmic problems that are natural
counterparts to the enumerative results:
\begin{itemize}
\item
How can one efficiently generate simple permutations in lexicographic
order?
\item
Is it possible to generate simple permutations uniformly at random in
worst-case linear time per permutation?
\item
How efficiently can one recognise a simple permutation?
\end{itemize}

With regards to the final question, there is a natural
dynamic programming algorithm that achieves the task in $O(n^2)$ time; so
the issue is whether one can do better.

\bibliographystyle{plain}

\begin{thebibliography}{10}

\bibitem{AA02}
M.H. Albert and M.D. Atkinson,
\newblock Simple permutations, partial well-order, and enumeration.
\newblock In M.H. Albert ed., {\em Proceedings, Permutation Patterns
2003}, \url{www.cs.otago.ac.nz/trseries/oucs-2003-02.pdf}, 2003, pp. 5--9.

\bibitem{AS01}
M.D. Atkinson and T.~Stitt,
\newblock Restricted permutations and the wreath product.
\newblock {\em Discrete Math.}, {\bf 259} (2002), 19--36.

\bibitem{BR01}
Edward~A. Bender and L.~Bruce Richmond,
\newblock An asymptotic expansion for the coefficients of some power series.
  {II}. {L}agrange inversion.
\newblock {\em Discrete Math.}, {\bf 50} (1984), 135--141.

\bibitem{Co01}
Louis Comtet,
\newblock {\em Advanced combinatorics}.
\newblock D. Reidel, 1974.

\bibitem{FS01}
P.~Flajolet and R.~Sedgewick,
\newblock {\em Analytic Combinatorics---Symbolic Combinatorics},
\newblock Preprint published electronically at,
\url{http://algo.inria.fr/flajolet/Publications/books.html}, 2002.

\bibitem{GJ01}
I.~P. Goulden and D.~M. Jackson,
\newblock {\em Combinatorial enumeration}.
\newblock John Wiley \& Sons, 1983.


\bibitem{Ka01}
Irving Kaplansky,
\newblock The asymptotic distribution of runs of consecutive elements,
\newblock {\em Ann. Math. Statistics}, {\bf 16} (1945), 200--203.

\bibitem{Kl01}
Martin Klazar,
\newblock Non-P-recursiveness of numbers of matchings or linear chord diagrams
  with many crossings,
\newblock {\em Adv. Appl. Math.}, {\bf 30} (2003), 126--136.

\bibitem{kumm} {E. E. Kummer,}
\newblock {\"Uber die Erg\"anzungss\"atze zu
den allgemeinen Reciprocit\"atsgezetzen},
\newblock {\em J. Reine
Angew. Math.}, {\bf 44} (1852), 93--146.

\bibitem{Mu01}
M.~M. Murphy,
\newblock {\em Restricted permutations, antichains, atomic classes and stack
  sorting},
\newblock PhD thesis, University of St. Andrews, 2002.

\bibitem{Si01}
David Singmaster,
\newblock Notes on binomial coefficients. {I}. {A} generalization of {L}ucas'
  congruence,
\newblock {\em J. London Math. Soc. (2)}, {\bf 8} (1974), 545--548.

\bibitem{S01}
N.J.A. Sloane,
\newblock The on-line encyclopedia of integer sequences,
\newblock \url{http://www.research.att.com/~njas/sequences/}, 2003.

\bibitem{St02}
Richard~P. Stanley,
\newblock {\em Enumerative combinatorics. {V}ol. 2}, volume~62 of {\em
  Cambridge Studies in Advanced Mathematics},
\newblock Cambridge University Press, 1999.


\end{thebibliography}

\end{document}